# A rational logit dynamic for decision-making under uncertainty: well-posedness, vanishing-noise limit, and numerical approximation

The first draft under review at some international conference

Hidekazu Yoshioka[1][0000-0002-5293-3246], Motoh Tsujimura[2][0000-0001-6975-9304], and Yumi Yoshioka[3][0000-0002-0855-699X]

[1] Japan Advanced Institute of Science and Technology, 1-1 Asahidai, Nomi, 923-1292, Japan
[2] Doshisha University, Karasuma-Higashi-iru, Imadegawa-dori, Kamigyo-ku, 602-8580, Japan
[3] Gifu University, Yanagido 1-1, Gifu, 501-1193, Japan
yoshih@jaist.ac.jp; mtsujimu@mail.doshisha.ac.jp;
yoshioka.yumi.k6@f.gifu-u.ac.jp

**Abstract.** The classical logit dynamic on a continuous action space for decision-making under uncertainty is generalized to the dynamic where the exponential function for the softmax part has been replaced by a rational one that includes the former as a special case. We call the new dynamic as the rational logit dynamic. The use of the rational logit function implies that the uncertainties have a longer tail than that assumed in the classical one. We show that the rational logit dynamic admits a unique measure-valued solution and the solution can be approximated using a finite difference discretization. We also show that the vanishing-noise limit of the rational logit dynamic exists and is different from the best-response one, demonstrating that influences of the uncertainty tail persist in the rational logit dynamic. We finally apply the rational logit dynamic to a unique fishing competition data that has been recently acquired by the authors.

## 1 Introduction

### 1.1 Research Background

Modeling social dynamics has been a central research topic to better understand and manage human interactions in the contemporary world. Social dynamics can be modelled through the dynamic game theory where interactions among (infinitely) many agents are described through an evolution equation. Analysis of this equation is therefore a central issue in the theory. Both trajectory and its equilibrium of a reasonable solution to the evolution equation are important to analyze transitions of the social state of interest. Such examples include but are not limited to the best-response dynamic [1], replicator dynamic [2], projection dynamic [3], gradient dynamic [4], and logit dynamic [5]. Evolutionary systems generalizing the above-mentioned dynamics have also been proposed along with their well-posedness and stability results [6,7,8].

We focus on the logit dynamic [5] in a continuous action space as it serves as the simplest evolution equation for decision-making of individuals in an uncertain environment. The logit dynamic reads the following evolution equation governing a time-dependent probability measure $\mu$ on a compact set $\Omega$ (the equation will be defined more rigorously in **Section 2**):

$$\frac{\mathrm{d}\mu}{\mathrm{d}t} = \frac{\exp(\eta^{-1}U(x;\mu))}{\int_\Omega \exp(\eta^{-1}U(y;\mu))\mathrm{d}y} - \mu \text{ for time } t > 0 \tag{1}$$

subject to a prescribed initial condition. Here, $U$ is a bounded utility depending on the social state $\mu$ as a probability measure on the collection of individuals' actions parameterized in $\Omega$, and $\eta > 0$ is a noise intensity. Phenomenologically, a larger $\eta$ implies a larger uncertainty that individuals face in the given environment during the decision-making. The noise acts to regularize the dynamic through the logit function. An equilibrium of (1) is interpreted as a maximizing $\mu$ of the utility $U$ under uncertainty.

The logit dynamic (1) is an evolution equation nonlocal in space. The first term in the right-hand side of (1) uses the logit function, which is here the exponential function so that the term serves as a softmax function. A stationary state of the dynamic (1) under the vanishing-noise limit $\eta \to +0$, if it exists, approximates a Nash equilibrium of the maximization problem of the utility $U$ with respect to probability measures $\mu$ [5]. The dynamic for a nonzero $\eta$ serves as a model of dynamic decision-making in an uncertain environment, which is of a separate interest.

It has recently been found that the use of a different logit function in the dynamic (1), such as a polynomial function, instead of the exponential function substantially affects the equilibrium profile of the dynamic. Lahkar et al. [9] numerically demonstrated that different logit functions lead to different equilibria for small $\eta$. Subsequently, Yoshioka [10] theoretically showed along with examples that the growth speed of the logit function controls the equilibrium; a different growth speed may give a different equilibrium. Their analysis clarified influences of the logit functions on the equilibria, whereas the analysis of time-dependent solutions is still scarce.

### 1.2   Objectives and Contributions

The objectives as well as contributions of this study are formulation, analysis, discretization, and computational application of a generalized logit dynamic: an evolution equation using a wider class of logit function than the classical one (1). As a logit function, we use the $\kappa$-exponential function with $\kappa \in [0,1]$, which is a positive, increasing, and strictly convex function [11]:

$$e_\kappa(z) = \left(\kappa z + \sqrt{\kappa^2 z^2 + 1}\right)^{\frac{1}{\kappa}}, \ z \in \mathbb{R} \tag{2}$$

with the convention $e_0(z) = \exp(z)$. The $\kappa$-exponential function is therefore a generalization of the exponential one. An elementary calculation yields

$$\frac{d e_\kappa(z)}{dz} = \frac{1}{\sqrt{\kappa^2 z^2 + 1}} e_\kappa(z), \ z \in \mathbb{R}, \qquad (3)$$

showing that the derivative of $e_\kappa$ is bounded in each compact set of $\mathbb{R}$.

We demonstrate that the resulting logit dynamic, which we call the rational logit dynamic because it uses a rational(-like) function (2) in the softmax part, admits a unique measure-valued solution. Moreover, we explicitly obtain the evolution equation corresponding to the vanishing-noise limit $\eta \to +0$ along with a convergence result, with which we can track solution trajectories with and without uncertainties. Rational logit functions also recently found their applications in accelerating the convergence in machine learning algorithms [12,13].

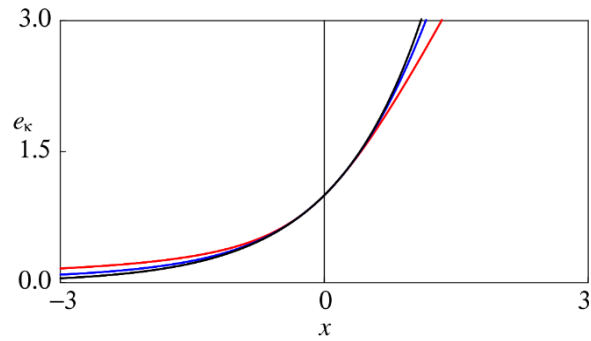

**Fig. 1.** Profiles of $\kappa$-exponential functions: $\kappa$ equals 0 (black), 0.5 (blue), and 1 (red).

We also propose a provably-convergence finite difference method for the rational logit dynamic. Finally, we apply the generalized logit dynamic to a unique fishing competition data that has been recently acquired by the authors, so that this study covers an engineering application.

The rest of this paper is organized as follows. **Section 2** formulates the rational logit dynamic, performs its mathematical analysis, and presents its finite difference discretization. **Section 3** shows a computational application of the rational logit dynamic. **Section 4** concludes this study and its future perspectives. **Appendix** contains a technical proof and the data used in our application.

## 2     Formulation and Analysis

### 2.1   Rational Logit Dynamic

We consider a decision-making problem of individuals parameterized by the closed unit interval $\Omega = [0,1]$. Borel algebra on $\Omega$ is denoted as $\mathfrak{B}$. The set of probability measures on $\Omega$ is denoted as $\mathfrak{M}$. The variational norm $\|\cdot\|$ for $\mu \in \mathfrak{M}$ is given by [5]

$$\|\mu\| = \sup_g \left| \int_\Omega g(y) \mu(\mathrm{d}y) \right|, \ g : \Omega \to [-1,1] \text{ is measurable.} \quad (4)$$

The difference $\|\mu - \nu\|$ for $\mu, \nu \in \mathfrak{M}$ is defined as in (4) where $\mu$ is formally replaced by $\mu - \nu$. The utility $U : \Omega \times \mathfrak{M} \to \mathbb{R}$ is bounded and Lipschitz continuous:

$$|U(x;\mu) - U(y;\nu)| \leq K(|x-y| + \|\mu - \nu\|) \text{ for any } x, y \in \Omega \text{ and } \mu, \nu \in \mathfrak{M} \quad (5)$$

with a constant $K > 0$ not depending on $x, y, \mu, \nu$. More technically rigorously, we should consider a utility $U : \Omega \times \mathfrak{M} \to \mathbb{R}$, where $\mathfrak{F}$ ($\mathfrak{M} \subset \mathfrak{F}$) is the set of finite signed measures on $\Omega$. However, this is possible for our and several important cases (see (15)) by assuming that the constant $K$ in (5) depends on $\max\{\|\mu\|, \|\nu\|\}$.

The rational logit dynamic in this paper is the evolution equation governing the time-dependent probability measure $\mu : [0, +\infty) \times \mathfrak{B} \to \mathfrak{M}$:

$$\frac{\mathrm{d}\mu(t,A)}{\mathrm{d}t} = \frac{\int_A e_\kappa(\eta^{-1} U(y;\mu)) \mathrm{d}y}{\int_\Omega e_\kappa(\eta^{-1} U(y;\mu)) \mathrm{d}y} - \mu(t,A) \text{ for any } A \in \mathfrak{B} \text{ and } t > 0 \quad (6)$$

subject to a prescribed initial condition $\mu(0, \cdot) \in \mathfrak{M}$.

We close this subsection with several remarks. The time derivative in (6) is understood in the strong sense with the variational norm (e.g., Mendoza-Palacios and [2]). The Lipschitz condition (5) is a typical assumption for logit dynamics [9,10], and is actually not so restrictive as the constant $K$ is allowed to be arbitrary large if it is bounded. The $\kappa$-exponential function seems to be less famous than the $q$-exponential one often used in the Tsallis formalism [14]:

$$e_{(q)}(z) = (1 + (1-q)z)^{\frac{1}{1-q}}, \ z \in \mathbb{R} \text{ if } 1 + (1-q)z > 0 \quad (7)$$

with a parameter $q > 0$. The $q$-exponential function reduces to the exponential one at $q = 1$ like the $\kappa$-exponential one. The $q$-exponential function is explicitly connected to heavy-tailed noises [10] particularly the generalized extreme value one [15], and therefore admits a firm background in the context of logit dynamics; however, a serious disadvantage of the $q$-exponential function is that it is not defined over $\mathbb{R}$. The $\kappa$- exponential function completely avoids this theoretical issue. As shown later, what is important in the rational logit dynamic (6) is the growth speed of the logit function, namely the function $e_\kappa\left(\eta^{-1}(\cdot)\right)$, for a large argument: $e_\kappa(z) \sim O\left(z^{1/\kappa}\right)$ for a sufficiently large $z > 0$. In this view, a similar result will be expected to follow when using the $q$-exponential function with the choice $q = 1 - \kappa \in (0,1)$. The dynamic (6) still governs a maximizer of the utility $U$ under uncertainty where the noise has a heavier tail than that in the classical one (1), such that the heaviness increases as $\kappa$ does.

### 2.2  Well-posedness and Stability

The unique existence of solutions to the rational logit dynamic follows by a direct application of Proof of Theorem 3.4 in Larkar and Riedel [5]. Here, a solution to a (rational) logit dynamic is a time-dependent probability measure $\mu(t,\cdot) \in \mathfrak{M}$ for $t \in [0, +\infty)$ that is continuously differentiable with respect to $t > 0$. The only difference between their and our logit dynamics are the logit functions. To apply their result, it suffices to show the following technical inequality for any $x \in \Omega$ and $\mu, \nu \in \mathfrak{M}$:

$$\left| \frac{e_\kappa\left(\eta^{-1}U(x;\mu)\right)}{\int_\Omega e_\kappa\left(\eta^{-1}U(y;\mu)\right)\mathrm{d}y} - \frac{e_\kappa\left(\eta^{-1}U(x;\nu)\right)}{\int_\Omega e_\kappa\left(\eta^{-1}U(y;\nu)\right)\mathrm{d}y} \right| \leq L\|\mu - \nu\|. \tag{8}$$

Here, $L > 0$ is a constant independent from $x \in \Omega$ and $\mu, \nu \in \mathfrak{M}$. The inequality (8) follows from the Lipschitz continuity $U$ assumed in (5) and its boundedness, the Lipschitz continuity of the $\kappa$-exponential function in each compact set in $\mathbb{R}$ that follows from (3), strict positivity of the $\kappa$-exponential function in each compact set in $\mathbb{R}$ that follows from its definition, and the fact that a composition of functions that are Lipschitz continuous is again a Lipschitz continuous function. With the help of these continuity conditions, we can use the generalized Picard–Lindelöf theorem (e.g., Theorem A.3 in Lahkar et al. [9]). Consequently, we obtain **Proposition 1** below.

### *Proposition 1*

*The rational logit dynamic (6) given a prescribed initial condition $\mu(0,\cdot) \in \mathfrak{M}$ admits a unique solution $\mu(t,\cdot) \in \mathfrak{M}$ for $t \in [0, +\infty)$.*

Having proven the unique existence of solutions to the rational logit dynamic, we analyze its vanishing-noise limit $\eta \to +0$. With the notations analogous to (6), the limit equation, the dynamic obtained under this limit, is inferred by formally taking the limit $\eta \to +0$ in the logit function considering (2):

$$\frac{\mathrm{d}\mu(t,A)}{\mathrm{d}t} = \frac{\int_A \max\{U(y;\mu),0\}^{1/\kappa}\mathrm{d}y}{\int_\Omega \max\{U(y;\mu),0\}^{1/\kappa}\mathrm{d}y} - \mu(t,A) \text{ for any } A \in \mathfrak{B} \text{ and } t > 0 \tag{9}$$

with the initial condition being unchanged. The following proposition is the main theoretical contribution of this study, which states a stability of the rational logit dynamic (6) with respect to $\eta > 0$ as well as its convergence to the limit one (9).

*Proposition 2*
*Assume that $U(x;\mu) \geq \underline{U}$ with a constant $\underline{U} > 0$ for any $x \in \Omega$ and $\mu \in \mathfrak{M}$. Fix constants $\eta_0 > 0$ and $T > 0$. Solutions to (6) with $\eta \in (0, \eta_0]$ and (9) with the same prescribed initial condition are denoted as $\mu_\eta$ and $\mu_0$, respectively. Then, it follows that*

$$\lim_{\eta \to +0} \sup_{0 \leq t \leq T} \|\mu_0(t,\cdot) - \mu_\eta(t,\cdot)\| = 0. \qquad (10)$$

See, **Appendix** for the proof of **Proposition 2**. Note that the constants $\eta_0, T > 0$ can be taken to be arbitrary large if they are bounded. **Proposition 2** therefore shows that the limit equation (9) indeed serves as a vanishing-noise limit of (6). The positivity assumption of $U$ was imposed for a technical reason, but it may be removed in some case as computationally demonstrated in **Section 3**. We can also obtain the parameter continuity result below:

$$\lim_{\eta_2 \to \eta_1} \sup_{0 \leq t \leq T} \|\mu_{\eta_1}(t,\cdot) - \mu_{\eta_2}(t,\cdot)\| = 0 \qquad (11)$$

with $\eta \in (0, \eta_0]$ under the assumption of **Proposition 2**. To show this, it suffices to see the triangle inequality

$$\|\mu_{\eta_1}(t,\cdot) - \mu_{\eta_2}(t,\cdot)\| \leq \|\mu_{\eta_1}(t,\cdot) - \mu_0(t,\cdot)\| + \|\mu_0(t,\cdot) - \mu_{\eta_2}(t,\cdot)\|. \qquad (12)$$

### 2.3   Numerical Discretization

We present a finite difference discretization of the rational logit dynamic (6). Set the resolution $N \in \mathbb{N}$ with $N \geq 2$. The domain $\Omega$ is uniformly discretized into $N$ cells

$$\Omega_i = \left[(i-1)/N, i/N\right) \ (i = 1, 2, ..., N-1) \text{ and } \Omega_N = [1 - 1/N, 1]. \qquad (13)$$

The midpoint of $\Omega_i$ is denoted as $x_{i-1/2}$. The dynamic (6) with the choice $A = \Omega_i$ is discretized as follows:

$$\frac{d\hat{\mu}(t,\Omega_i)}{dt} = \frac{N^{-1} e_\kappa\left(\eta^{-1} U(x_{i-1/2}; \hat{\mu})\right)}{\sum_{j=1}^{N} N^{-1} e_\kappa\left(\eta^{-1} U(x_{j-1/2}; \hat{\mu})\right)} - \hat{\mu}(t,\Omega_i) \text{ for } i = 1, 2, ..., N \text{ and } t > 0. \qquad (14)$$

We used the fact that the length of each cell is $N^{-1}$ and $\hat{\mu}(t,\Omega_i)$ is the discretization of $\mu(t,\Omega_i)$. The initial condition $\mu(0,\cdot)$ will also be discretized in each cell.

To complete the spatial discretization (14), we need to specify the form of the utility $U$. For that purpose, we assume the following widely-used form [2,9,10]:

$$U(x;\mu) = \int_\Omega f(x,y)\mu(\mathrm{d}y) \text{ for any } x \in \Omega \text{ and } \mu \in \mathfrak{M} \qquad (15)$$

with a function $f : \Omega \times \Omega \to \mathbb{R}$ that is continuous on the compact set $\Omega \times \Omega$, and hence bounded as well as uniformly continuous on $\Omega \times \Omega$. We apply the discretization

$$U(x;\mu) = \int_\Omega f(x,y)\mu(\mathrm{d}y) \to \sum_{k=1}^N f(x_{j-1/2}, x_{k-1/2})\hat{\mu}(t,\Omega_k) = U(x_{j-1/2};\hat{\mu}) \qquad (16)$$

in (14). The full discretization of the equation (14) is finally obtained by applying the common first-order forward Euler method that evaluates the right-hand side of (14) explicitly. The limit equation (9) can also be discretized using the same finite difference method by directly utilizing (15)-(16). By using each $\hat{\mu}(t,\Omega_i)$, the space semi-discretized solution to the rational logit dynamic is potentially defined as the time-dependent measure $\mu_N$ satisfying the equality:

$$\mu_N(t,A) = \int_A N \sum_{i=1}^N \mathbb{I}_{\Omega_i}(y)\hat{\mu}(t,\Omega_i)\mathrm{d}y \text{ for any } A \in \mathfrak{B} \text{ and } t \in [0,+\infty), \qquad (17)$$

where $\mathbb{I}_S(\cdot)$ is the indicator function of set $S$ such that $\mathbb{I}_S(y) = 1$ if $y \in S$ and is 0 otherwise. This $\mu_N$ is checked to be a probability measure at each $t \in [0,+\infty)$ by substituting $A = \Omega$ in (17) considering the identity $\int_{\Omega_i} N\mathrm{d}y = 1$.

We state that the solution to the semi-discretized dynamic (14) exists and converges to a solution to the original one. The proof of the proposition below is omitted here because it is based on Proof of Proposition 3.1 of Lahkar et al. [16], where it suffices to use the Lipschitz continuity like that for **Proposition 1**.

*Proposition 3*
*Assume $\kappa(0,1]$ and (15) with $f : \Omega \times \Omega \to \mathbb{R}$ that is continuous on $\Omega \times \Omega$. Solutions to (6) and (14) with initial conditions $\mu(0,\cdot), \mu_N(0,\cdot)$ are denoted as $\mu$ and $\mu_N$, respectively. Also assume that $\lim_{N \to +\infty} \|\mu_N(0,\cdot) - \mu(0,\cdot)\| = 0$. Then, it follows that*

$$\lim_{N \to +\infty} \|\mu_N(t,\cdot) - \mu(t,\cdot)\| = 0 \text{ at each } t > 0. \qquad (18)$$

Note that a similar result applies to the limit equation. Hence, the assumption $\eta > 0$ can be formally replaced by $\eta \geq 0$ in **Proposition 3**. Finally, this proposition still holds true if a sufficiently regular term satisfying the Lipschitz continuity of the form (5) is added to the right-hand side of (15) (see, (19)-(20) in the next section).

## 3     Application

We apply the rational logit dynamic (6) and the finite difference method to an engineering problem related to inland fisheries. The authors have been studying environmental and ecological dynamics of the Hii River, San-in area, Japan since 2015. This river is a major habitat of the migratory fish *Plecoglossus altivelis altivelis*, which is one of the most famous inland fishery resources in Japan (e.g., Murase and Iguchi [17]). Toami (i.e., casting net) is a major fishing method of the fish *P. altivelis*.

The Hii River Fishery Cooperative (HRFC) authorizing inland fisheries in the Hii River is holding the Toami competition in each year. Participants of the competition are resisted in pairs. The total number of pairs in each competition is 10 to 16 in recent years. The registered pairs compete the total catch of the fish *P. altivelis* by the Toami angling in the two hours (typically 10:00 a.m. to 12:00 a.m.) in a fixed day. All participants must use casting nets under a common regulation. Top few pairs (about three pairs in each year) are awarded by the HRFC. **Table 2** in **Appendix** summaries the number of catches of the fish *P. altivelis* in each pair in each Toami competition.

The total number of pairs and their catches are different among different competitions. We therefore normalize the data of the fish catches in each competition by the division by the maximum (e.g., 82 in 2023 and 43 in 2016). The fish catches in each competition then belong to $\Omega = [0,1]$, to which our formulation applies. In this view, the fish catch can be understood as the normalized harvesting effort. The data used in our application is an ensemble of the normalized data in each year in **Table 2**.

We model the Toami competition by the following utility

$$U(x;\mu) = \int_\Omega \left\{ -ax^2 + b|x-y|^c \right\} \mu(\mathrm{d}y) + d \max\left\{ \alpha - \int_x^1 \mu(\mathrm{d}y), 0 \right\} \quad (19)$$

with $a,b,c,d \geq 0$ and $\alpha \in (0,1)$. Each term in the right-hand side of (19) is explained as follows. The first term represents the harvesting cost that is assumed to be quadratic with respect to the effort. The quadratic assumption of the cost is typical in the control and optimization problems. The second term represents the individual difference modelled in a way that this term is larger when the fish catches between different pairs are more different, and hence is expected to induce a competition. The third term represents the awarding scheme by HRFC for the Toami competition that the $100\alpha\%$ upper tier is awarded. This term is continuous in $x$ if $\mu$ has a density. The utility (19) models a competition with an award mechanism. Rigorously, to apply our mathematical framework to (19), the last term should be regularized. We use the following continuous regularization of the step function:

$$\int_x^1 \mu(\mathrm{d}y) = \int_0^1 \mathbb{I}_{(x,1]}(y) \mu(\mathrm{d}y) \to \int_0^1 \max\left\{ 0, \min\left\{ 1, \frac{y-x+\varepsilon}{\varepsilon} \right\} \right\} \mu(\mathrm{d}y) \text{ with } \varepsilon > 0. \quad (20)$$

Indeed, we have the Lipschitz continuity

$$\left| \int_0^1 \max\left\{0, \min\left\{1, \frac{y-x+\varepsilon}{\varepsilon}\right\}\right\} (\mu(\mathrm{d}y) - \nu(\mathrm{d}y)) \right| \leq \|\mu - \nu\| + \varepsilon^{-1} |x - z| \qquad (21)$$

for any $x, z \in \Omega$ and $\mu, \nu \in \mathfrak{M}$. In the sequel, we choose $\varepsilon = N^{-1}$.

We need to identify the model parameters in (19) as well as those in the rational logit dynamic. We set $d = 1$ without significant loss of generality. Further, we set $c = 1$ to simplify the problem. We choose $\alpha = 0.2$ considering the award scheme employed by the HRFC. The remaining parameters $a, b, \eta, \kappa$ are determined using a trial-and-error approach so that the sum of the squares of the relative errors of the average and standard deviation of $x$ is minimized between the empirical and computational results. The computational resolution to discretize the rational logit dynamic and its limit equation is $N = 500$ with the time step $\Delta t = 0.001$. The initial condition is the uniform distribution in $\Omega$. The solution to the discretized dynamic is judged to be stationary when the computed probability density function (PDF) given by $N\hat{\mu}$ satisfies

$$\max_{1 \leq i \leq N} \left| n\hat{\mu}(t_{k+1}, \Omega_i) - n\hat{\mu}(t_k, \Omega_i) \right| \leq \delta \left(= 10^{-11}\right) \text{ at some } k \geq 0. \qquad (22)$$

The stationary solution is then given by $\hat{\mu}\left(t_{k+1}, \Omega_{(\cdot)}\right)$ with the smallest $k$ satisfying the condition (22). We assume the stationary state against the empirical data.

The identified parameter values are $a = 0.27$, $b = 0.23$, $\eta = 0.01$, $\kappa = 1.0$. The empirical average and standard deviation of the normalized harvesting effort $x$ are 0.32471 and 0.30352, respectively; the theoretical average and standard deviation using the fitted model are 0.32471 and 0.30377, respectively. The relative errors of the average and standard deviation between the empirical and theoretical ones are 0.0000011 and 0.00025, respectively, that are sufficiently small. **Fig. 2** compares the empirical and theoretical PDFs of the effort, suggesting their reasonable agreement. Particularly, the model captures the bimodal nature of the empirical data.

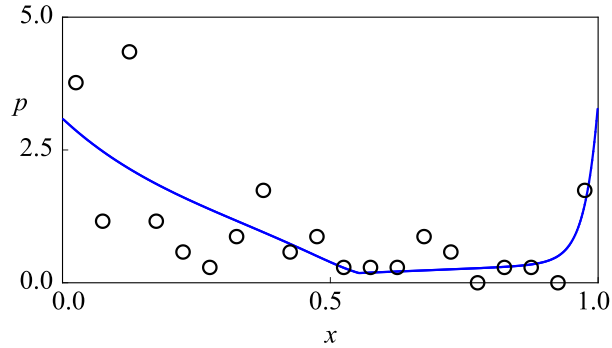

**Fig. 2.** Comparison of the empirical PDF (black) and theoretical one (blue).

Having identified the model from the data, we conduct a sensitivity analysis of the rational logit dynamic. We focus on the two parameters characterizing the dynamics, which are the noise intensity $\eta$ and shape parameter $\kappa$.

Firstly, we analyze $\eta$. **Fig. 3** shows the computed profiles of the PDF $p$ for different values of $\eta$ including the extreme case $\eta = 0$, where the other parameters remain the same with those identified from the data. As indicated in **Fig. 3**, increasing $\eta$ flattens the profile of $p$, which is in accordance with the fact that $\eta$ is the intensity of noises blurring the information available for anglers. The bimodal nature of $p$ is preserved during increasing $\eta$. The profiles between the case $\eta = 0$ and $\eta = 0.01$ are visually almost the same with each other. We then evaluate the convergence speed of the numerical solutions with $\eta > 0$ to that of $\eta = 0$ using the maximum norm at time $t = 1$, $t = 10$, and $t = 100$ as shown in **Table 1**. The convergence speed is slightly better than the first order in $\eta$, not violating the theory (32).

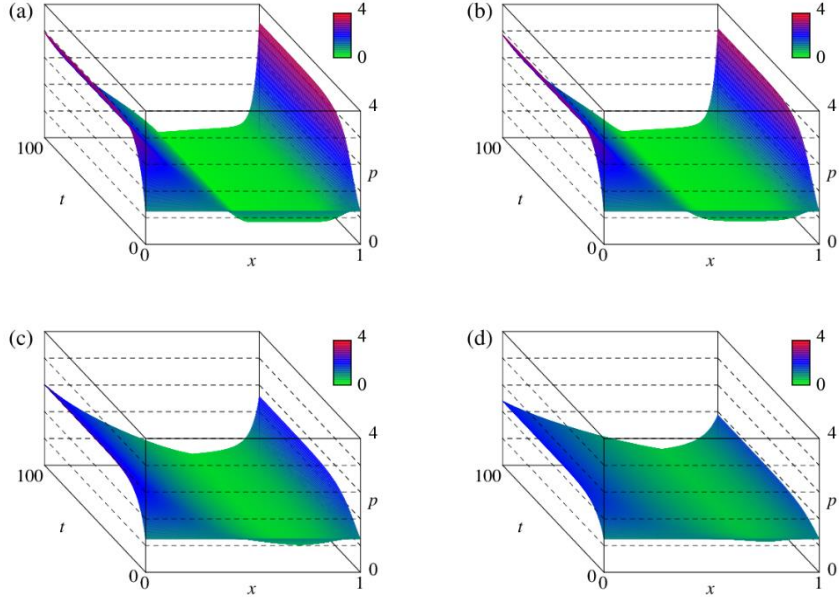

**Fig. 3.** The transient profiles of the PDF: (a) $\eta = 0$ (b) $\eta = 0.01$, (c) $\eta = 0.05$, (d) $\eta = 0.1$.

**Table 1.** The maximum norm error between the solutions with $\eta > 0$ and that with $\eta = 0$.

|  |  | $\eta$ |  |  |  |
|---|---|---|---|---|---|
|  |  | 0.0001 | 0.001 | 0.01 | 0.1 |
| $t = 1$ | Error | 8.33E-05 | 3.74E-03 | 1.15E-01 | 1.03E+00 |
|  | Convergence rate | 1.65.E+00 | 1.49.E+00 | 9.50.E-01 |  |
| $t = 10$ | Error | 2.44E-05 | 2.40E-03 | 1.73E-01 | 1.94E+00 |
|  | Convergence rate | 1.99.E+00 | 1.86.E+00 | 1.05.E+00 |  |

Secondly, we analyze influences of the shape parameter $\kappa$ on solutions to the rational logit dynamic. **Fig. 4** compares stationary PDFs of the logit dynamic for different values of $\kappa$, where the other parameters remain the same with those identified from the data. The maximum values of $p$ in the right-most cell are 14.60 for $\kappa = 0.1$ and 16.01 for $\kappa = 0$, which are placed above the figure panel. Decreasing $\kappa$ from the identified value 1 toward 0 sharpens the profile of stationary $p$, while maintaining the flat part that exists around $0.6 \leq x \leq 0.9$. The stationary PDFs are entirely positive for both large and small $\kappa$, and hence the support of $p$ is identical to $\Omega$. Studying a more complex problem with a support expansion/contraction will be an interesting problem.

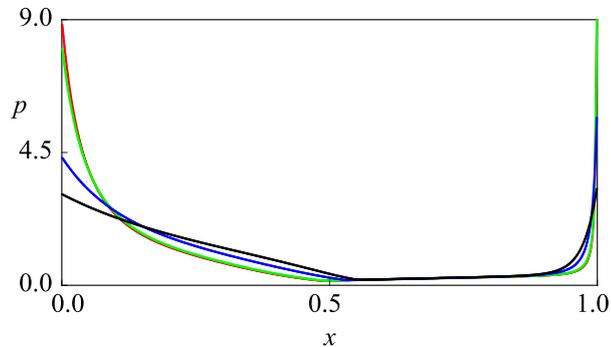

**Fig. 4.** Comparison of the stationary PDFs: $\kappa = 0$ (red), $\kappa = 0.1$ (green), $\kappa = 0.5$ (blue), $\kappa = 1$ (black).

## 4    Conclusion

We proposed a rational logit dynamic where a rational function based on the $\kappa$-exponential function replaced the exponential function in the soft max part. It was proven that the rational logit dynamic is well-posed and has a parameter continuity with respect to the noise intensity that may vanish. We pointed out that the vanishing-noise limit of the rational logit dynamic depends on the shape parameter $\kappa$ of the logit function, while such a phenomenon is not encountered in the classical logit dynamic, thereby the rational logit dynamic should be understood as a qualitatively different model from the classical one. We finally demonstrated a unique application of the rational logit dynamic to the fishing competition data.

The proposed logit function can be generalized because our analysis implied that what determines behavior of the rational logit dynamic under a small-noise limit is not the whole profile of the logit function but its asymptotic for a large argument. We will generalize the proposed mathematical framework based on this conjecture. The proposed logit function grows slower than the exponential one, while the case where it grows faster than the exponential one, like $\exp(x^2)$, will also be worth investigating.

Examining the logit dynamics against real problems are of importance to deeper study their applicability and limitations. In particular, the regularity condition assumed the utility, such as the Lipschitz continuity with respect to the action and probability measure, should also be investigated in future so that the logit dynamic can be formulated under a weaker assumption. Establishment of a provably convergent as well as more efficient numerical method to discretize the logit dynamics will be a critical issue as well. Exploring the connection of the logit dynamics to mean filed games [18] has been under investigations by the authors.

**Appendix**
*Proof of Proposition 2*
The key estimate in the proof of **Proposition 2** is the following; given $\kappa \in (0,1]$, $\mu \in \mathfrak{M}$, $\eta \in (0, \eta_0]$ with a constant $\eta_0 > 0$, by the boundedness of $U$ it follows that

$$\left| (\eta/(2\kappa))^{1/\kappa} e_\kappa \left( \eta^{-1} U(x;\mu) \right) - U(x;\mu)^{1/\kappa} \right|$$
$$= \left| \left( U(x;\mu)/2 + \sqrt{U(x;\mu)^2/4 + \eta^2 \kappa^{-2}/4} \right)^{\frac{1}{\kappa}} - U(x;\mu)^{1/\kappa} \right| \leq C\eta, \quad (23)$$

where $C > 0$ **is a constant independent from** $\eta$. This independence is crucial in our context, particularly when considering the limit $\eta \to +0$ (see (27)-(28) below).

Now, for any $t \in (0,T]$ and $A \in \mathfrak{B}$, it follows that

$$\frac{\mathrm{d}}{\mathrm{d}t} \left( \mu_0(t,A) - \mu_\eta(t,A) \right)$$
$$= \frac{\int_A U(y;\mu_0)^{1/\kappa} \mathrm{d}y}{\int_\Omega U(y;\mu_0)^{1/\kappa} \mathrm{d}y} - \frac{(\eta/(2\kappa))^{1/\kappa} \int_A e_\kappa \left( \eta^{-1} U(y;\mu_\eta) \right) \mathrm{d}y}{(\eta/(2\kappa))^{1/\kappa} \int_\Omega e_\kappa \left( \eta^{-1} U(y;\mu_\eta) \right) \mathrm{d}y} - \left( \mu_0(t,A) - \mu_\eta(t,A) \right). \quad (24)$$

We have the following estimate with a constant $C_1 > 0$ independent from $\eta, \mu_0, \mu_\eta$:

$$\int_\Omega U(y;\mu_0)^{1/\kappa} \mathrm{d}y (\eta/(2\kappa))^{1/\kappa} \int_\Omega e_\kappa \left( \eta^{-1} U(y;\mu_\eta) \right) \mathrm{d}y$$
$$\geq \int_\Omega U(y;\mu_0)^{1/\kappa} \mathrm{d}y \int_\Omega U(y;\mu_\eta)^{1/\kappa} \mathrm{d}y > C_1 > 0 \quad (25)$$

We also have the following estimate at each $y \in \Omega$:

$$(\eta/(2\kappa))^{1/\kappa} \int_\Omega e_\kappa\left(\eta^{-1} U(x;\mu_\eta)\right) \mathrm{d}x U(y;\mu_0)^{1/\kappa} - \int_\Omega U(x;\mu_0)^{1/\kappa} \mathrm{d}x (\eta/(2\kappa))^{1/\kappa} e_\kappa\left(\eta^{-1} U(y;\mu_\eta)\right)$$

$$= (\eta/(2\kappa))^{1/\kappa} \int_\Omega e_\kappa\left(\eta^{-1} U(x;\mu_\eta)\right) \mathrm{d}x U(y;\mu_0)^{1/\kappa} - \int_\Omega U(x;\mu_0)^{1/\kappa} \mathrm{d}x U(y;\mu_0)^{1/\kappa}$$

$$+ \int_\Omega U(x;\mu_0)^{1/\kappa} \mathrm{d}x U(y;\mu_0)^{1/\kappa} - \int_\Omega U(x;\mu_0)^{1/\kappa} \mathrm{d}x (\eta/(2\kappa))^{1/\kappa} e_\kappa\left(\eta^{-1} U(y;\mu_\eta)\right)$$

$$= \int_\Omega \left((\eta/(2\kappa))^{1/\kappa} e_\kappa\left(\eta^{-1} U(x;\mu_\eta)\right) - U(x;\mu_0)^{1/\kappa}\right) \mathrm{d}x U(y;\mu_0)^{1/\kappa}$$

$$+ \int_\Omega U(x;\mu_0)^{1/\kappa} \mathrm{d}x \left\{U(y;\mu_\eta) - (\eta/(2\kappa))^{1/\kappa} e_\kappa\left(\eta^{-1} U(y;\mu_\eta)\right)\right\}$$

$$\leq C_2 \left\{ \begin{array}{l} \int_\Omega \left|(\eta/(2\kappa))^{1/\kappa} e_\kappa\left(\eta^{-1} U(x;\mu_\eta)\right) - U(x;\mu_0)^{1/\kappa}\right| \mathrm{d}x \\ + \left|U(y;\mu_\eta) - (\eta/(2\kappa))^{1/\kappa} e_\kappa\left(\eta^{-1} U(y;\mu_\eta)\right)\right| \end{array} \right\}$$

$$\leq C_2 \left\{ \begin{array}{l} (\eta/(2\kappa))^{1/\kappa} \int_\Omega \left|e_\kappa\left(\eta^{-1} U(x;\mu_\eta)\right) - e_\kappa\left(\eta^{-1} U(x;\mu_0)\right)\right| \mathrm{d}x \\ + \int_\Omega \left|(\eta/(2\kappa))^{1/\kappa} e_\kappa\left(\eta^{-1} U(x;\mu_0)\right) - U(x;\mu_0)^{1/\kappa}\right| \mathrm{d}x \\ + \left|U(y;\mu_\eta) - (\eta/(2\kappa))^{1/\kappa} e_\kappa\left(\eta^{-1} U(y;\mu_\eta)\right)\right| \end{array} \right\}$$

(26)

with a constant $C_2 > 0$ independent from $\mu_0, \mu_\eta$. By (23) and $U(x;\mu_\eta) > 0$, in (26) we obtain

$$\left|U(y;\mu_\eta) - (\eta/(2\kappa))^{1/\kappa} e_\kappa\left(\eta^{-1} U(y;\mu_\eta)\right)\right| \leq C\eta, \quad (27)$$

$$\int_\Omega \left|(\eta/(2\kappa))^{1/\kappa} e_\kappa\left(\eta^{-1} U(x;\mu_0)\right) - U(x;\mu_0)^{1/\kappa}\right| \mathrm{d}x \leq C\eta, \quad (28)$$

and

$$(\eta/(2\kappa))^{1/\kappa} \int_\Omega \left|e_\kappa\left(\eta^{-1} U(x;\mu_\eta)\right) - e_\kappa\left(\eta^{-1} U(x;\mu_0)\right)\right| \mathrm{d}x \leq C_3(\eta_0) \|\mu_\eta - \mu_0\| \quad (29)$$

with a constant $C_3(\eta_0) > 0$ independent from $\eta, \mu_0, \mu_\eta$. By (24)-(29), we obtain

$$\frac{\mathrm{d}}{\mathrm{d}t}\left(\mu_0(t,A) - \mu_\eta(t,A)\right) \leq C_4(\eta_0)\left(\eta + \|\mu_\eta(t,\cdot) - \mu_0(t,\cdot)\|\right) \quad (30)$$

with a constant $C_4(\eta_0) > 0$ independent from $\eta, \mu_0, \mu_\eta$. Hence, by integrating (30) for $(0,t)$, and taking the variational norm yields

$$\|\mu_\eta(t,\cdot) - \mu_0(t,\cdot)\| \leq \int_0^t C_4(\eta_0)\left(\eta + \|\mu_\eta(s,\cdot) - \mu_0(s,\cdot)\|\right) \mathrm{d}s. \quad (31)$$

Applying a classical Gronwall lemma to (31) yields

$$\left\| \mu_\eta(t,\cdot) - \mu_0(t,\cdot) \right\| \leq C_5(T,\eta_0)\eta \exp(C_5(\eta_0)T) \quad (32)$$

with a constant $C_5(T,\eta_0) > 0$ depending on $\eta_0, T$ but not on $\mu_0, \mu_\eta$. The conclusion (10) directly follows from (32).

□

*Collected data used in Section 3*

The number of catches of the fish *P. altivelis* in each pair in each Toami competition is summarized in the ascending order in **Table 2**. We consider that this kind of fish catch data is useful because it can be utilized not only for our study but also for other studies by other researchers. Members of each pair were anonymized. The Toami competition has been basically held by HRFC in each summer. It was not held in 2020, 2021, 2022 due to the outbreak of the coronavirus disease 2019. The data before 2015 may exist but was not available for us.

**Table 2.** The number of catches of the fish *P. altivelis* in each group in each year.

| 2016 Aug 7 | 2017 Aug 6 | 2018 Aug 5 | 2019 Aug 4 | 2023 Jul 30 |
|---|---|---|---|---|
| 1 | 0 | 2 | 0 | 0 |
| 1 | 5 | 6 | 1 | 0 |
| 2 | 5 | 6 | 1 | 3 |
| 2 | 8 | 8 | 2 | 9 |
| 4 | 9 | 17 | 4 | 9 |
| 5 | 16 | 17 | 4 | 9 |
| 6 | 17 | 21 | 4 | 10 |
| 6 | 21 | 23 | 5 | 12 |
| 8 | 22 | 24 | 6 | 12 |
| 10 | 25 | 53 | 16 | 13 |
| 14 | 28 | | 25 | 24 |
| 16 | 41 | | 28 | 30 |
| 21 | 42 | | 29 | 31 |
| 32 | | | 35 | 56 |
| 36 | | | 41 | 82 |
| 43 | | | | |

**Acknowledgments.** The authors would like to express their gratitude towards the officers and members of HRFC for their supports on our field surveys. This study was supported by JSPS grants No. 22K14441 and No. 22H02456.

**Disclosure of Interests.** The authors have no competing interests to disclose.